\documentclass[12pt]{article}
\usepackage{amsmath,amsfonts,amssymb,amscd}
\usepackage[enableskew]{youngtab}
\usepackage{theorem}
\usepackage[usenames]{color}
\nonstopmode

\newtheorem{theorem}{Theorem}[section]
\newtheorem{proposition}[theorem]{Proposition}

\newtheorem{definition}[theorem]{Definition}
\newtheorem{corollary}[theorem]{Corollary}
\newtheorem{conjecture}[theorem]{Conjecture}

\theorembodyfont{\rmfamily}
\newtheorem{example}[theorem]{Example}

\begin{document}
$\,$\vspace{10mm}

\begin{center}
{\textsf {\Huge Generalized Energy Statistics and}}
\vspace{3mm}\\
{\textsf {\Huge Kostka--Macdonald Polynomials}}
\vspace{15mm}\\
{\textsf{\LARGE  ${}^{\mbox{\small a}}$Anatol N. Kirillov and
${}^{\mbox{\small b}}$Reiho Sakamoto}}
\vspace{20mm}\\
{\textsf {${}^{\mbox{\small{a}}}$Research Institute for Mathematical Sciences,}}
\vspace{-1mm}\\
{\textsf {Kyoto University, Sakyo-ku, Kyoto, 606-8502, Japan}}
\vspace{-1mm}\\
{\textsf {kirillov@kurims.kyoto-u.ac.jp}}
\vspace{3mm}\\
{\textsf {${}^{\mbox{\small{b}}}$Department of Physics, Tokyo University of Science,}}
\vspace{-1mm}\\
{\textsf {Kagurazaka, Shinjuku-ku, Tokyo, 162-8601, Japan}}
\vspace{-1mm}\\
{\textsf {reiho@rs.kagu.tus.ac.jp}}
\vspace{30mm}
\end{center}

\begin{quotation}
\noindent
{\bf Abstract:}
We give an interpretation of the $t=1$ specialization of the
modified Macdonald polynomial as a generating function of the
energy statistics defined on the
set of paths arising in the context of Box-Ball Systems
(BBS-paths for short).
We also introduce one parameter generalizations of the
energy statistics on the set of BBS-paths which all,
conjecturally, have the same distribution.
\bigskip\\
{\bf R\'{e}sum\'{e}:}
Nous donnons une int\'{e}rpr\'{e}tation de la sp\'{e}cialisation 
\`{a} $t=1$ du polyn\^{o}me de Macdonald modifi\'{e} comme fonction
g\'{e}n\'{e}ratrice des statistiques d'\'{e}nergie d\'{e}finies sur
l'ensemble des chemins qui apparaissent dans la th\'{e}orie
des Syst\`{e}mes BBS (BBS-chemins).
Nous pr\'{e}sentons \'{e}galement des g\'{e}n\'{e}ralisations
\`{a} un param\`{e}tre de la statistique d'\'{e}nergie sur les chemins BBS qui toutes,
conjecturalement, ont la m\^{e}me distribution.\bigskip
\end{quotation}

\pagebreak

\section{Introduction}
The purpose of the present paper is two-fold.
First of all we would like to draw attention to a rich
combinatorics hidden behind the dynamics of Box-Ball Systems,
and secondly, to connect the former with the theory of
modified Macdonald polynomials.
More specifically, our final goal is to give an interpretation
of the Kostka--Macdonald polynomials $K_{\lambda,\mu}(q,t)$
as a {\it refined partition function} of a certain box-ball systems
depending on initial data $\lambda$ and $\mu$.

Box-Ball Systems (BBS for short) were invented by
Takahashi--Satsuma \cite{TS,Tak} as a wide class of discrete
integrable soliton systems.
In the simplest case, BBS are described by simple combinatorial
procedures using boxes and balls.
One can see the simplest but still very interesting examples
of the BBS by the free software available at \cite{SakWol}.
Despite its simple outlook, it is known that
the BBS have various remarkably deep properties:
\begin{itemize}
\item
Local time evolution rule of the BBS coincides with the isomorphism
of the crystal bases \cite{HHIKTT,FOY}.
Thus the BBS possesses quantum integrability.
\item
BBS are ultradiscrete (or tropical) limit of the
usual soliton systems \cite{TTMS,KSY}.
Thus the BBS possesses classical integrability
at the same time.
\item
Inverse scattering formalism of the BBS \cite{KOSTY} coincides with
the rigged configuration bijection originating in
completeness problem of the Bethe states \cite{Kir1,KR}, see also \cite{Sak3}.
\end{itemize}
Let us say a few words about the main results of this note.
\begin{itemize}
\item
We will identify the space of states of a  BBS with the corresponding weight subspace in the 
tensor product of fundamental (or rectangular) representations of the Lie algebra 
$\mathfrak{gl}(n).$
\item
In the case of statistics {\it tau},
our main result can be formulated as a computation
of the corresponding  partition function for the BBS in terms of the
values of the Kostka--Macdonald polynomials at $t=1.$
\item
In the case of the statistics {\it energy},
our result can be formulated as 
an interpretation of the corresponding
partition function for the BBS as the $q$-weight 
multiplicity of a certain irreducible representation of the Lie 
algebra $\mathfrak{gl}(n)$ in the tensor product of
the fundamental representations.  We {\it expect}
that the same statement is valid for the BBS 
corresponding to the tensor product of rectangular representations.

\hspace{3mm}Let us remind that a {\it $q$-analogue of the multiplicity} of a highest 
weight $\lambda$ in the tensor product $\bigotimes_{a=1}^{L}~V_{s_{a} \omega_{r_a}}$
of the highest weight $s_a~\omega_{r_a},$~$a=1,\ldots,L,$ irreducible representations 
$V_{s_{a} \omega_{r_{a}}}$ of the Lie algebra $\mathfrak {gl}(n)$  is defined as
$$q\mbox{-Mult }[V_{\lambda}~ :~ \bigotimes_{a=1}^{L}~V_{s_{a} \omega_{r_{a}}}] =
\sum_{\eta}~K_{\eta, R}~K_{\eta, \lambda}(q),$$
where  $K_{\eta, R}$ stands for the parabolic Kostka number corresponding to the 
sequence of rectangles $R:= \{ (s_{a}^{r_{a}}) \}_{a=1,\ldots,L},$ 
see e.g. \cite{Kir}, \cite{KiSh}.
\end{itemize} 

A combinatorial description of the modified Macdonald polynomials
has been obtained by Haglund--Haiman--Loehr \cite{HHL}.
In Section \ref{sec:haglund} we give an interpretation of two Haglund's
statistics in the context of the box-ball systems, i.e.,
in terms of the BBS-paths.
Namely, we identify the set of BBS paths of weight $\alpha$
with the set $\mathcal{P}(\alpha)$ which is the weight $\alpha$
component in the tensor product of crystals corresponding
to vector representations.
We have observed that from the proof given in \cite{HHL}
one can prove the following identity
\begin{equation}\label{eq:intro}
\sum_{p\in\mathcal{P}(\alpha)}
q^{\mathrm{inv}_\mu (p)}
t^{\mathrm{maj}_\mu (p)}=
\sum_{\eta\vdash |\mu|}
K_{\eta,\alpha}
\tilde{K}_{\eta,\mu}(q,t),
\end{equation}
see Proposition \ref{prop:reducedHHL} and Corollary \ref{cor:reducedHHL}.
One of the main problems we are interested in is to generalize
the identity Eq.(\ref{eq:intro}) on more wider set of the BBS-paths.

Our result about connections of the energy partition
functions for BBS and $q$-weight 
multiplicities suggests a deep hidden connections
between partition functions for the BBS 
and characters of the Demazure modules, solutions to the
$q$-difference Toda equations, cf.\cite{GLO}, ... . 

As an interesting open problem we want to give raise
a question about an interpretation 
of the sums $\sum_{\eta}~K_{\eta,R}~K_{\eta,\lambda}(q,t),$ where
 $K_{\eta,\lambda}(q,t)$ 
denotes the Kostka--Macdonald polynomials \cite{M},
as {\it refined partition functions} for the BBS 
corresponding to the tensor product of rectangular representations 
$R=\{(s_a^{r_a}) \}_{1 \le a \le n}$.
In other words, one can ask: what is a meaning of the 
second statistics (see \cite{HHL}) in the Kashiwara theory \cite{Kas}  of crystal 
bases (of type A) ?

This paper is abbreviated and updated version of our paper \cite{KiSa}.
The main novelty of the present paper is the definition of a one parameter
family of statistics on the set of BBS-paths which generalizes those introduced
in \cite{KiSa}, see Conjecture \ref{con:main}.
It conjecturally gives a new family of MacMahonian
statistics on the set of transportation matrices, see \cite{Kir}.

Organization of the present paper is as follows.
In Section \ref{sec:KRcrystal} we outlook the basic
definitions and facts related to the Kashiwara's theory of crystal base
in the case of type $A_n^{(1)}$.
We also remind definitions of the combinatorial $R$-matrix and
definition of the energy function.
We illustrate definitions by simple example.
In Section \ref{sec:enestat}, we introduce the energy statistics
and the set of the BBS.
In Section \ref{sec:BBS} we remind definition of box-ball systems
and state some of their simplest properties.
In Section \ref{sec:haglund} we remind definition of the Haglund's
statistics and give their interpretation in terms of the BBS-paths.
Sections \ref{sec:HHL} and \ref{sec:gentau} contain our main results
and conjectures.
In particular it is not difficult to see that Haglund's
statistics maj${}_\mu$ and inv${}_\mu$ do not compatible with
the Kostka--Macdonald polynomials for general partitions
$\lambda$ and $\mu$.
In Section \ref{sec:HHL} we state a conjecture which describes
the all pairs of partitions $(\lambda,\mu)$ for those
the restriction of the Haglund--Haiman--Loehr formula on
the set of highest weight paths of shape $\mu$
coincide with the Kostka--Macdonald polynomial $\tilde{K}_{\lambda,\mu}(q,t)$.

\section{Kirillov--Reshetikhin crystal}\label{sec:KRcrystal}
\subsection{$A^{(1)}_n$ type crystal}
Let $W_s^{(r)}$ be a $U_q'(\mathfrak{g})$ Kirillov--Reshetikhin
module, where we shall consider the case $\mathfrak{g}=A_n^{(1)}$.
The module $W_s^{(r)}$ is indexed by a Dynkin node
$r\in I=\{1,2,\ldots,n\}$ and $s\in\mathbb{Z}_{>0}$.
As a $U_q(A_n)$-module, $W_s^{(r)}$ is isomorphic to the irreducible
module corresponding to the partition $(s^r)$.
For arbitrary $r$ and $s$, the module $W_s^{(r)}$
is known to have crystal bases \cite{Kas,KMN2}, 
which we denote by $B^{r,s}$.
As the set, $B^{r,s}$ is consisting of all column strict
semi-standard Young tableaux
of depth $r$ and width $s$ over the alphabet $\{1,2,\ldots,n+1\}$.

For the algebra $A_n$, let $P$ be the weight lattice,
$\{\Lambda_i\in P|i\in I\}$ be the fundamental roots,
$\{\alpha_i\in P|i\in I\}$ be the simple roots,
and $\{h_i\in \mathrm{Hom}_\mathbb{Z}(P,\mathbb{Z})|i\in I\}$
be the simple coroots.
As a type $A_n$ crystal, $B=B^{r,s}$ is equipped with
the Kashiwara operators $\tilde{e}_i,\tilde{f}_i:B\longrightarrow B\cup\{0\}$
and $\mathrm{wt}:B\longrightarrow P$ ($i\in I$) satisfying
\begin{eqnarray*}
&&\tilde{f}_i(b)=b'\Longleftrightarrow
  \tilde{e}_i(b')=b\quad \mbox{ if }b,b'\in B,\\
&&\mathrm{wt}\big(\tilde{f}_i(b)\big)=\mathrm{wt}(b)-\alpha_i
\quad \mbox{ if }\tilde{f}_i(b)\in B,\\
&&\langle h_i,\mathrm{wt}(b)\rangle =\varphi_i(b)-\varepsilon_i(b).
\end{eqnarray*}
Here $\langle\cdot ,\cdot\rangle$ is the natural pairing
and we set 
$\varepsilon_i(b)=\max\{m\ge0\mid \tilde{e}_i^m b\ne0\}$
and $\varphi_i(b)=\max\{m\ge0\mid \tilde{f}_i^m b\ne0\}$.
Actions of the Kashiwara operators
$\tilde{e}_i$, $\tilde{f}_i$
for $i\in I$ coincide with the one described in \cite{KN}.
Since we do not use explicit forms of these operators,
we omit the details.
See \cite{O} for complements of this section.
Note that in our case $A_n$, we have $P=\mathbb{Z}^{n+1}$
and $\alpha_i=\epsilon_i-\epsilon_{i+1}$ where
$\epsilon_i$ is the $i$-th canonical unit vector of $\mathbb{Z}^{n+1}$.
We also remark that $\mathrm{wt}(b)=(\lambda_1,\cdots,\lambda_{n+1})$
is the weight of $b$, i.e., $\lambda_i$ counts the number of letters
$i$ contained in tableau $b$.

For two crystals $B$ and $B'$,
one can define the tensor product
$B\otimes B'=\{b\otimes b'\mid b\in B,b'\in B'\}$.
The actions of the Kashiwara operators on tensor
product have simple form.
Namely, the operators 
$\tilde{e}_i,\tilde{f}_i$ act on $B\otimes B'$ by
\begin{eqnarray*}
\tilde{e}_i(b\otimes b')&=&\left\{
\begin{array}{ll}
\tilde{e}_i b\otimes b'&\mbox{ if }\varphi_i(b)\ge\varepsilon_i(b')\\
b\otimes \tilde{e}_i b'&\mbox{ if }\varphi_i(b) < \varepsilon_i(b'),
\end{array}\right. \\
\tilde{f}_i(b\otimes b')&=&\left\{
\begin{array}{ll}
\tilde{f}_i b\otimes b'&\mbox{ if }\varphi_i(b) > \varepsilon_i(b')\\
b\otimes \tilde{f}_i b'&\mbox{ if }\varphi_i(b)\le\varepsilon_i(b'),
\end{array}\right.
\end{eqnarray*}
and $\mathrm{wt}(b\otimes b')=\mathrm{wt}(b)+
\mathrm{wt}(b')$.
We assume that $0\otimes b'$ and $b\otimes 0$ as $0$.
Then it is known that there is the unique crystal isomorphism
$R:B^{r,s}\otimes B^{r',s'}
\stackrel{\sim}{\rightarrow}B^{r',s'}\otimes B^{r,s}$.
We call this map (classical) combinatorial $R$
and usually write the map $R$ simply by $\simeq$.

Let us consider the affinization of the crystal $B$.
As the set, it is
\begin{equation}
\mathrm{Aff}(B)=\{b[d]\, |\, b\in B,\, d\in\mathbb{Z}\}.
\end{equation}
There is also explicit algorithm for
actions of the affine Kashiwara operators $\tilde{e}_0$, $\tilde{f}_0$
in terms of  the promotion operator \cite{Shimo}.
For the tensor product
$b[d]\otimes b'[d']\in
\mathrm{Aff}(B)\otimes\mathrm{Aff}(B')$,
we can lift the (classical)
combinatorial $R$ to affine case
as follows:
\begin{equation}\label{eq:affineR}
b[d]\otimes b'[d']\stackrel{R}{\simeq}
\tilde{b}'[d'-H(b\otimes b')]\otimes
\tilde{b}[d+H(b\otimes b')],
\end{equation}
where $b\otimes b'\simeq \tilde{b}'\otimes \tilde{b}$
is the isomorphism of (classical) combinatorial $R$.
The function $H(b\otimes b')$ is called the energy function
and defined by a certain set of axioms.
We will give explicit forms of the combinatorial $R$ and
energy function in the next section.

\subsection{Combinatorial $R$ and energy function}
We give an explicit description of the combinatorial $R$-matrix
(combinatorial $R$ for short)
and energy function on $B^{r,s}\otimes B^{r',s'}$.
To begin with we define few terminologies about Young tableaux.
Denote rows of a Young tableaux $Y$ by $y_1,y_2,\ldots y_r$
from top to bottom.
Then row word $row(Y)$ is defined by concatenating rows as
$row(Y)=y_ry_{r-1}\ldots y_1$.
Let $x=(x_1,x_2,\ldots )$ and $y=(y_1,y_2,\ldots )$ be two partitions.
We define concatenation of $x$ and $y$ by the partition
$(x_1+y_1,x_2+y_2,\ldots )$.

\begin{proposition}[\cite{Shimo}]\label{prop:shimozono}
$b\otimes b'\in B^{r,s}\otimes B^{r',s'}$ is mapped to
$\tilde{b}'\otimes \tilde{b}\in B^{r',s'}\otimes B^{r,s}$
under the combinatorial $R$, i.e.,
\begin{equation}
b\otimes b'\stackrel{R}{\simeq}\tilde{b}'\otimes\tilde{b},
\end{equation}
if and only if
\begin{equation}
(b'\leftarrow row(b))=(\tilde{b}\leftarrow row(\tilde{b}')).
\end{equation}
Moreover, the energy function $H(b\otimes b')$ is given by
the number of nodes of $(b'\leftarrow row(b))$
outside the concatenation of partitions
$(s^r)$ and $({{s'}^{r'}})$.
\end{proposition}

For special cases of $B^{1,s}\otimes B^{1,s'}$,
the function $H$ is called unwinding number
in \cite{NY}.
Explicit values for the case
$b\otimes b'\in B^{1,1}\otimes B^{1,1}$
are given by $H(b\otimes b')=
\chi (b<b')$
where $\chi(\mathrm{True})=1$ and
$\chi(\mathrm{False})=0$.

In order to describe the
algorithm for finding $\tilde{b}$ and $\tilde{b}'$ from
the data $(b'\leftarrow row(b))$,
we introduce a terminology.
Let $Y$ be a tableau, and $Y'$ be a subset of $Y$
such that $Y'$ is also a tableau.
Consider the set theoretic subtraction $\theta =Y\setminus Y'$.
If the number of nodes contained in $\theta$ is $r$
and if the number of nodes of $\theta$ contained in
each row is always 0 or 1,
then $\theta$ is called vertical $r$-strip.

Given a tableau
$Y=(b'\leftarrow row(b))$, let $Y'$ be the upper left
part of $Y$ whose shape is $(s^r)$.
We assign numbers from 1 to $r's'$
for each node contained in $\theta =Y\setminus Y'$
by the following procedure.
Let $\theta_1$ be the vertical $r'$-strip of $\theta$
as upper as possible.
For each node in $\theta_1$,
we assign numbers 1 through $r'$ from the bottom to top.
Next we consider $\theta\setminus\theta_1$,
and find the vertical $r'$ strip $\theta_2$
by the same way.
Continue this procedure until all nodes of $\theta$
are assigned numbers up to $r's'$.
Then we apply inverse bumping procedure according
to the labeling of nodes in $\theta$.
Denote by $u_1$ the integer which is ejected
when we apply inverse bumping procedure starting
from the node with label 1.
Denote by $Y_1$ the tableau such that
$(Y_1\leftarrow u_1)=Y$.
Next we apply inverse bumping procedure
starting from the node of $Y_1$ labeled by 2,
and obtain the integer $u_2$ and tableau $Y_2$.
We do this procedure until we obtain $u_{r's'}$
and $Y_{r's'}$.
Finally, we have
\begin{equation}
\tilde{b}'=(\emptyset\leftarrow u_{r's'}u_{r's'-1}\cdots
u_1),\qquad
\tilde{b}=Y_{r's'}.
\end{equation}\begin{example}
Consider the following tensor product:
$$b\otimes b'=
\Yvcentermath1
\young(114,236)\otimes\young(23,34,45)\in
B^{2,3}\otimes B^{3,2}.$$
{}From $b$, we have $row(b)=236114$, hence we have
$$\left
(
\Yvcentermath1\young(23,34,45)
\leftarrow 236114
\right)=
\newcommand{\yonsan}{4_3}
\newcommand{\sanni}{3_2}
\newcommand{\yonichi}{4_1}
\newcommand{\goyon}{5_4}
\newcommand{\yongo}{4_5}
\newcommand{\sanroku}{3_6}
\Yvcentermath1
\young(113\yonsan ,226,\sanroku\sanni ,\yongo\yonichi ,\goyon) .
$$
Here subscripts of each node indicate the order of
inverse bumping procedure.
For example, we start from the node $4_1$ and obtain
$$\left
(
\Yvcentermath1
\young(1234,236,34,4,5)
\leftarrow 1
\right)=
\Yvcentermath1
\young(1134,226,33,44,5),
\qquad\mathrm{therefore,}\qquad
Y_1=
\newcommand{\yonsan}{4_3}
\newcommand{\yonni}{4_2}
\newcommand{\goyon}{5_4}
\newcommand{\yongo}{4_5}
\newcommand{\sanroku}{3_6}
\Yvcentermath1
\young(123\yonsan ,236,\sanroku\yonni ,\yongo ,\goyon)
,\qquad
u_1=1.
$$
Next we start from the node $4_2$ of $Y_1$.
Continuing in this way, we obtain
$u_6u_5\cdots u_1=321421$ and
$Y_6=\Yvcentermath1
\young(334,456)$.
Since $(\emptyset\leftarrow 321421)=
\Yvcentermath1\young(11,22,34)$,
we obtain
$$\Yvcentermath1
\young(114,236)\otimes\young(23,34,45)\simeq
\young(11,22,34)\otimes\young(334,456)\, ,
\qquad
H\left
(
\Yvcentermath1
\young(114,236)\otimes\young(23,34,45)
\right)
=3.
$$
Note that the energy function is derived from the
concatenation of shapes of $b$ and $b'$,
i.e., $\Yvcentermath1\yng(5,5,2)\,$.
\end{example}

\section{Energy statistics and its generalizations on the set of paths}
\label{sec:enestat}
For a path $b_1\otimes b_2\otimes
\cdots\otimes b_L\in B^{r_1,s_1}\otimes
B^{r_2,s_2}\otimes\cdots\otimes B^{r_L,s_L}$,
let us define elements
$b_j^{(i)}\in B^{r_j,s_j}$ for $i<j$
by the following isomorphisms of the combinatorial $R$;
\begin{eqnarray}
&&b_1\otimes b_2\otimes\cdots\otimes
b_{i-1}\otimes b_{i}\otimes\cdots\otimes
b_{j-1}\otimes b_j\otimes\cdots\nonumber\\
&\simeq&
b_1\otimes b_2\otimes\cdots\otimes
b_{i-1}\otimes b_{i}\otimes\cdots\otimes
b_j^{(j-1)}\otimes b_{j-1}'\otimes\cdots
\nonumber\\
&\simeq&\cdots\nonumber\\
&\simeq&
b_1\otimes b_2\otimes\cdots\otimes
b_{i-1}\otimes b_{j}^{(i)}\otimes
\cdots\otimes
b_{j-2}'\otimes b_{j-1}'\otimes\cdots,
\end{eqnarray}
where we have written
$b_k\otimes b_j^{(k+1)}\simeq
b_j^{(k)}\otimes b_k'$
assuming that $b_j^{(j)}=b_j$.

Define the statistics $\mathrm{maj}(p)$ by
\begin{equation}\label{def:maj}
\mathrm{maj}(p)=
\sum_{i<j}H(b_i\otimes b_j^{(i+1)}).
\end{equation}
For example, consider a path
$a=a_1\otimes a_2\otimes\cdots\otimes a_L
\in (B^{1,1})^{\otimes L}$.
In this case, we have $a_j^{(i)}=a_i$,
since the combinatorial $R$ act on
$B^{1,1}\otimes B^{1,1}$ as identity.
Therefore, we have
\begin{equation}
\mathrm{maj}(a)
=\sum_{i=1}^{L-1}(L-i)\chi (a_i<a_{i+1}).
\end{equation}

Define another statistics {\it tau} as follows.
\begin{definition}\label{def:generaltau}
For the path $p\in B^{r_1,s_1}\otimes
B^{r_2,s_2}\otimes\cdots\otimes B^{r_L,s_L}$,
define $\tau^{r,s}$ by
\begin{equation}
\tau^{r,s}(p)=\mathrm{maj}
(u^{(r)}_s\otimes p),
\end{equation}
where $u^{(r)}_s$ is the highest element of $B^{r,s}$.
\end{definition}
Here the highest element
$u_{s}^{(r)}\in B^{r,s}$ is the tableau
whose $i$-th row is occupied by integers $i$.
For example, $u_4^{(3)}=\Yvcentermath1
\young(1111,2222,3333)\,$.
In particular, the statistics
$\tau^{r,1}$ on $B^{1,1}$ type paths
$a\in (B^{1,1})^{\otimes L}$ has the following form;
\begin{equation}\label{eq:taumu}
\tau^{r,1} (a)=L\cdot\chi (r<a_1)+
\sum_{i=1}^{L-1}(L-i)\chi (a_i<a_{i+1}),
\end{equation}
where $a_1$ denotes the first letter of the path $a$.
Note that $\tau^{1,1}$ is a special case of the tau functions
for the box-ball systems \cite{KSY,Sak1}
which originates as an ultradiscrete limit of
the tau functions for the KP hierarchy \cite{JM}.
\begin{definition}\label{def:taumu}
For composition $\mu=(\mu_1,\mu_2,\cdots,\mu_n)$,
write $\mu_{[i]}=\sum_{j=1}^i \mu_j$
with convention $\mu_{[0]}=0$.
Then we define a generalization
of $\tau^{r,1}$ by
\begin{equation}
\tau_\mu^{r,1}(a)=\sum_{i=1}^{n}
\tau^{r,1}(a_{[i]}),
\end{equation}
where
\begin{equation}
a_{[i]}=
a_{\mu_{[i-1]}+1}\otimes
a_{\mu_{[i-1]}+2}\otimes\cdots\otimes
a_{\mu_{[i]}}
\in (B^{1,1})^{\otimes\mu_i}.
\end{equation}
\end{definition}
Note that we have
$a=a_{[1]}\otimes a_{[2]}\otimes\cdots\otimes
a_{[n]}$, i.e., the path $a$ is partitioned
according to $\mu$.

\section{Box-ball system}\label{sec:BBS}
In this section, we summarize basic facts about
the box-ball system in order to explain physical
origin of $\tau^{1,1}$.
For our purpose, it is convenient to
express the isomorphism of the combinatorial
$R$: $a\otimes b\simeq b'\otimes a'$
by the following vertex diagram:
\begin{center}
\unitlength 13pt
\begin{picture}(4,4)
\put(0,2.0){\line(1,0){3.2}}
\put(1.6,1.0){\line(0,1){2}}
\put(-0.6,1.8){$a$}
\put(1.4,0){$b'$}
\put(1.4,3.2){$b$}
\put(3.4,1.8){$a'$}
\put(4.1,1.7){.}
\end{picture}
\end{center}
Successive applications of the combinatorial $R$
is depicted by concatenating these vertices.

Following \cite{HHIKTT,FOY}, we define time evolution
of the box-ball system $T^{(a)}_l$.
Let $u_{l,0}^{(a)}=u_{l}^{(a)}\in B^{a,l}$ be the highest
element and $b_i\in B^{r_i,s_i}$.
Define $u_{l,j}^{(a)}$ and $b_i'\in B^{r_i,s_i}$
by the following diagram.
\begin{equation}\label{def:T_l}
\unitlength 13pt
\begin{picture}(22,5)(0,-0.5)
\multiput(0,0)(5.8,0){2}{
\put(0,2.0){\line(1,0){4}}
\put(2,0){\line(0,1){4}}
}
\put(-1.4,1.8){$u_{l,0}^{(a)}$}
\put(1.7,4.2){$b_1$}
\put(1.7,-0.8){$b_1'$}
\put(4.2,1.8){$u_{l,1}^{(a)}$}
\put(7.5,4.2){$b_2$}
\put(7.5,-0.8){$b_2'$}
\put(10.0,1.8){$u_{l,2}^{(a)}$}
\multiput(11.5,1.8)(0.3,0){10}{$\cdot$}
\put(14.7,1.8){$u_{l,L-1}^{(a)}$}
\put(17,0){
\put(0,2.0){\line(1,0){4}}
\put(2,0){\line(0,1){4}}
}
\put(18.7,4.2){$b_L$}
\put(18.7,-0.8){$b_L'$}
\put(21.2,1.8){$u_{l,L}^{(a)}$}
\end{picture}
\end{equation}
$u_{l,j}^{(a)}$ are usually called {\it carrier}
and we set $u_{l,0}^{(a)}:=u_{l}^{(a)}$.
Then we define operator $T_l^{(a)}$ by
\begin{equation}
T_l^{(a)}(b)=b'=
b_1'\otimes b_2'\otimes\cdots\otimes b_L'.
\end{equation}
Recently \cite{Sak3}, operators $T^{(a)}_l$ have used
to derive crystal theoretical meaning of the
rigged configuration bijection.

It is known (\cite{KOSTY} Theorem 2.7) that there
exists some $l\in\mathbb{Z}_{>0}$ such that
\begin{equation}
T^{(a)}_{l}=T^{(a)}_{l+1}=T^{(a)}_{l+2}=\cdots
(=:T^{(a)}_\infty ).
\end{equation}
If the corresponding path is $b\in (B^{1,1})^{\otimes L}$,
we have the following combinatorial description
of the box-ball system \cite{TS,Tak}.
We regard $\fbox{1}\in B^{1,1}$ as an empty box of
capacity 1, and $\fbox{$i$}\in B^{1,1}$
as a ball of label (or internal degree of freedom) $i$
contained in the box.
Then we have:
 
\begin{proposition}[\cite{HHIKTT}]\label{prop:bbs}
For a path $b\in (B^{1,1})^{\otimes L}$ of type $A^{(1)}_n$,
$T^{(1)}_\infty (b)$ is given by
the following procedure.
\begin{enumerate}
\item
Move every ball only once.

\item
Move the leftmost ball with label $n+1$
to the nearest right empty box.

\item
Move the leftmost ball with label $n+1$
among the rest to its nearest right
empty box.

\item
Repeat this procedure until all of the balls
with label $n+1$ are moved.

\item
Do the same procedure 2--4 for the balls with
label $n$.

\item
Repeat this procedure successively until all
of the balls with label $2$ are moved.
\end{enumerate}
\end{proposition}
There are extensions of this box and ball algorithm
corresponding to generalizations of the box-ball
systems with respect to each affine Lie algebra,
see e.g., \cite{HKT3}.
Using this box and ball interpretation,
our statistics $\tau^{1,1} (b)$ admits the following
interpretation.

\begin{theorem}[\cite{KSY} Theorem 7.4]\label{th:tau=rho}
For a path $b\in (B^{1,1})^{\otimes L}$ of type $A^{(1)}_n$,
$\tau^{1,1} (b)$ coincides with number of all balls
$2,\cdots,n+1$ contained in paths
$b$, $T^{(1)}_\infty (b)$, $\cdots$,
$( T^{(1)}_\infty )^{L-1} (b)$.
\end{theorem}

\begin{example}
Consider the path $p=a\otimes b$ where
$a=4311211111$, $b=4321111111$.
Note that we omit all frames of tableaux of $B^{1,1}$
and symbols for tensor product.
We compute $\tau_{(10,10)}(p)$ by
using Theorem \ref{th:tau=rho}.
According to Proposition \ref{prop:bbs},
the time evolutions of the paths $a$ and
$b$ are as follows:
\begin{center}
\fbox{\hspace{-3pt}
$\begin{array}{llllllllll}
 4 & 3 & 1 & 1 & 2 & 1 & 1 & 1 & 1 & 1\\
 1 & 1 & 4 & 3 & 1 & 2 & 1 & 1 & 1 & 1\\
 1 & 1 & 1 & 1 & 4 & 1 & 3 & 2 & 1 & 1\\
 1 & 1 & 1 & 1 & 1 & 4 & 1 & 1 & 3 & 2\\
 1 & 1 & 1 & 1 & 1 & 1 & 4 & 1 & 1 & 1\\
 1 & 1 & 1 & 1 & 1 & 1 & 1 & 4 & 1 & 1\\
 1 & 1 & 1 & 1 & 1 & 1 & 1 & 1 & 4 & 1\\
 1 & 1 & 1 & 1 & 1 & 1 & 1 & 1 & 1 & 4
\end{array}$}\qquad\qquad
\fbox{\hspace{-3pt}
$\begin{array}{llllllllll}
 4 & 3 & 2 & 1 & 1 & 1 & 1 & 1 & 1 & 1\\
 1 & 1 & 1 & 4 & 3 & 2 & 1 & 1 & 1 & 1\\
 1 & 1 & 1 & 1 & 1 & 1 & 4 & 3 & 2 & 1\\
 1 & 1 & 1 & 1 & 1 & 1 & 1 & 1 & 1 & 4\\
 1 & 1 & 1 & 1 & 1 & 1 & 1 & 1 & 1 & 1\\
 1 & 1 & 1 & 1 & 1 & 1 & 1 & 1 & 1 & 1\\
 1 & 1 & 1 & 1 & 1 & 1 & 1 & 1 & 1 & 1\\
 1 & 1 & 1 & 1 & 1 & 1 & 1 & 1 & 1 & 1
\end{array}
$}
\end{center}
Here the left and right tables correspond to
$a$ and $b$, respectively.
Rows of left (resp. right) table represent
$a$, $T^{(1)}_\infty (a)$, $\cdots$,
$( T^{(1)}_\infty )^L (a)$ (resp., those for $b$)
from top to bottom.
Counting letters 2, 3 and 4 in each table,
we have $\tau^{1,1} (a)=16$, $\tau^{1,1} (b)=10$
and we get $\tau_{(10,10)}^{1,1}(p)=16+10=26$,
which coincides with the computation by Eq.(\ref{eq:taumu}).
Meanings of the above two dynamics
corresponding to paths $a$ and $b$ are
summarized as follows:
\begin{enumerate}
\item[$(a)$]
Dynamics of the path $a$.
In the first two rows, there are two solitons
(length two soliton $43$
and length one soliton 2),
and in the lower rows, there are also two solitons
(length one soliton 4 and length two soliton 32).
This is scattering of two solitons.
After the scattering, soliton 4 propagates
at velocity one and soliton 32 propagates
at velocity two without scattering.
\item[$(b)$]
Dynamics of the path $b$.
This shows free propagation of one soliton
of length three 432 at velocity three.
\end{enumerate}
\end{example}

\section{Haglund's statistics}\label{sec:haglund}
\paragraph{Tableaux language description}
For a given path
$a=a_1\otimes a_2\otimes\cdots\otimes a_L
\in (B^{1,1})^{\otimes L}$, associate tabloid $t$
of shape $\mu$
whose reading word coincides with $a$.
For example, to path
$p={\it abcdefgh}$
and the composition $\mu=(3,2,3)$
one associates the tabloid
\begin{equation}\label{eq:Haglundreading}
\Yvcentermath1
\young(cba,ed,hgf)\,\, .
\end{equation}
Denote the cell at the $i$-th row,
$j$-th column
(we denote the coordinate by $(i,j)$) of
the tabloid $t$ by $t_{ij}$.
Attacking region of the cell at $(i,j)$
is all cells
$(i,k)$ with $k<j$ or $(i+1,k)$ with $k>j$.
In the following diagram, gray zonal
regions are the attacking regions of
the cell $(i,j)$.

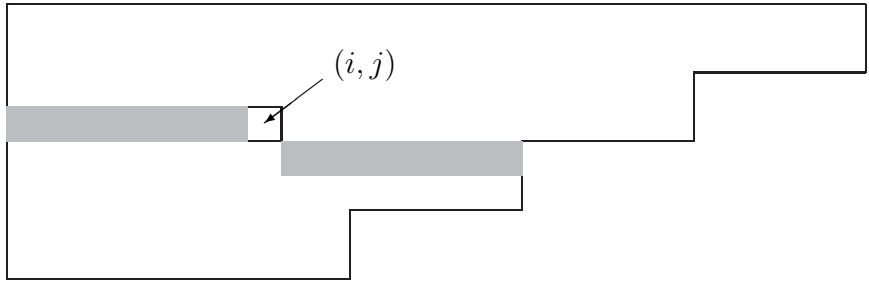
\unitlength 13pt
\begin{picture}(10,9)(-2,1.5)
\put(0,2){\line(0,1){8}}
\put(0,2){\line(1,0){5}}
\put(0,10){\line(1,0){25}}
\multiput(5,2)(5,2){4}{\line(1,0){5}}
\multiput(10,2)(5,2){4}{\line(0,1){2}}
\multiput(7,6)(0,1){2}{\line(1,0){1}}
\multiput(7,6)(1,0){2}{\line(0,1){1}}
\color[cmyk]{0,0,0,0.3}
\put(0,6){\rule{91pt}{13pt}}
\put(8,5){\rule{91pt}{13pt}}
\color{black}
\put(9.2,7.8){\vector(-4,-3){1.7}}
\put(9.5,8){$(i,j)$}
\end{picture}

\noindent
Follow \cite{HHL}, define $|\mathrm{Inv}_{ij}|$ by
\begin{equation}
|\mathrm{Inv}_{ij}|=
\#\{(k,l)\in \mbox{ attacking region for } (i,j)\, |\,
t_{kl}>t_{ij}\}.
\end{equation}
Then we define
\begin{equation}
|\mathrm{Inv}_\mu (a)|=\sum_{(i,j)\in \mu}
|\mathrm{Inv}_{ij}|.
\end{equation}

If we have $t_{(i-1)j}<t_{ij}$,
then the cell $(i,j)$ is called by {\it descent}.
Then define
\begin{equation}
{\rm Des}_\mu(a)=
\sum_{{\rm all\, descent\, }(i,j)}
(\mu_i-j).
\end{equation}
Note that $(\mu_i-j)$ is the arm length
of the cell $(i,j)$.

\paragraph{Path language description}
Consider two paths
$a^{(1)},a^{(2)}\in (B^{1,1})^{\otimes\mu}$.
We denote by
$a^{(1)}\otimes a^{(2)}=a_1\otimes a_2
\otimes\cdots\otimes a_{2\mu}$.
Then we define
\begin{equation}
\mathrm{Inv}_{(\mu ,\mu)}(a^{(1)},a^{(2)})=
\sum_{k=1}^{\mu}\sum_{i=k+1}^{k+\mu -1}
\chi (a_k<a_i).
\end{equation}
For more general cases
$a^{(1)}\in (B^{1,1})^{\otimes\mu_1}$ and
$a^{(2)}\in (B^{1,1})^{\otimes\mu_2}$
satisfying $\mu_1>\mu_2$, we define
\begin{equation}
\mathrm{Inv}_{(\mu_1,\mu_2)}(a^{(1)},a^{(2)}):=
\mathrm{Inv}_{(\mu_1,\mu_1)}
{(a^{(1)},1^{\otimes (\mu_1-\mu_2)}
\otimes a^{(2)})}.
\end{equation}
Then the above definition of $|\mathrm{Inv}_\mu(a)|$
is equivalent to
\begin{equation}
|\mathrm{Inv}_\mu(a)|=
\sum_{i=1}^{n-1}\mathrm{Inv}_{(\mu_i,\mu_{i+1})}.
\end{equation}

Consider two paths
$a^{(1)}\in (B^{1,1})^{\otimes\mu_1}$ and
$a^{(2)}\in (B^{1,1})^{\otimes\mu_2}$
satisfying $\mu_1\geq \mu_2$.
Denote $a=a^{(1)}\otimes a^{(2)}$.
Then define
\begin{equation}
\mathrm{Des}_{(\mu_1,\mu_2)}(a)=
\sum_{k=\mu_1-\mu_2+1}^{\mu_1}
(k-(\mu_1-\mu_2)-1)\chi (a_k<a_{k+\mu_2}).
\end{equation}
For the tableau $T$ of shape $\mu$ corresponding to
the path $a$, we define
\begin{equation}
\mathrm{Des}_\mu (T)=
\sum_{i=1}^n\mathrm{Des}_{(\mu_i,\mu_{i+1})}
(a_{[i]}\otimes a_{[i+1]}).
\end{equation}

\begin{definition}[\cite{Hag}]
For a path $a$,
statistics $\mathrm{maj}_\mu$ is defined by
\begin{equation}
\mathrm{maj}_\mu (a)=
\sum_{i=1}^{\mu_1}\mathrm{maj}(t_{1,i}\otimes t_{2,i}
\otimes\cdots\otimes t_{\mu'_i,i}).
\end{equation}
and $\mathrm{inv}_{\mu}(a)$ is defined by
\begin{equation}
\mathrm{inv}_\mu (a)=|\mathrm{Inv}_\mu(a)|-
\mathrm{Des}_\mu (a).
\end{equation}
\end{definition}

If we associate to a given path $p \in {\cal P}(\lambda)$
with the shape $\mu$ tabloid $T$, we sometimes write
$\mathrm{maj}_\mu (p)=\mathrm{maj}(T)$ and
$\mathrm{inv}_\mu (p)=\mathrm{inv}(T)$.

\section{Haglund--Haiman--Loehr formula}\label{sec:HHL}
Let $\tilde{H}_\mu (x;q,t)$ be the (integral form) modified Macdonald
polynomials where $x$ stands for
infinitely many variables $x_1,x_2,\cdots$.
Here $\tilde{H}_\mu (x;q,t)$ is obtained by
simple plethystic substitution (see, e.g.,
section 2 of \cite{Hai}) from the original
definition of the Macdonald polynomials \cite{M}.
Schur function expansion of $\tilde{H}_\mu (x;q,t)$
is given by
\begin{equation}
\tilde{H}_\mu (x;q,t)=
\sum_\lambda\tilde{K}_{\lambda ,\mu}(q,t)
s_\lambda (x),
\end{equation}
where $\tilde{K}_{\lambda ,\mu}(q,t)$
stands for the following transformation of the
Kostka--Macdonald polynomials:
\begin{equation}
\tilde{K}_{\lambda ,\mu}(q,t)=
t^{n(\mu)}K_{\lambda ,\mu}(q,t^{-1}).
\end{equation}
Here we have used notation
$n(\mu)=\sum_i(i-1)\mu_i$.
Then the celebrated
Haglund--Haiman--Loehr (HHL) formula is
as follows.

\begin{theorem}[\cite{HHL}]
Let $\sigma:\mu\rightarrow\mathbb{Z}_{>0}$
be the filling of the Young diagram $\mu$
by positive integers $\mathbb{Z}_{>0}$,
and define $x^\sigma=\prod_{u\in\mu}x_{\sigma(u)}$.
Then the Macdonald polynomial
$\tilde{H}_\mu (x;q,t)$ have the following
explicit formula:
\begin{equation}
\tilde{H}_\mu (x;q,t)=
\sum_{\sigma:\mu\rightarrow\mathbb{Z}_{>0}}
q^{\mathrm{inv}(\sigma)}
t^{\mathrm{maj}(\sigma)} x^\sigma .
\end{equation}
\end{theorem}

{}From the HHL formula,
we can show the following formula.
\begin{proposition}\label{prop:reducedHHL}
For any partition $\mu$ and composition $\alpha$ of the same size,
one has
\begin{equation}\label{eq:HHL}
\sum_{p\in\mathcal{P}(\alpha)}
q^{\mathrm{inv}_\mu (p)}
t^{\mathrm{maj}_\mu (p)}=
\sum_{\eta\vdash |\mu|}
K_{\eta,\alpha}
\tilde{K}_{\eta,\mu}(q,t),
\end{equation}
where $\mathcal{P}(\alpha )$ stands for the set of type $B^{1,1}$
paths of weight $\alpha=(\alpha_1,\alpha_2,\ldots,\alpha_{n+1})$
and $\eta$ runs over all partitions of size $|\mu|$.
\end{proposition}

\begin{corollary}\label{cor:reducedHHL}
The (modified) Macdonald polynomial $\tilde{H}_\mu (x;q,t)$
have the following
expansion in terms of the monomial symmetric functions
$m_{\lambda}(x)$:
\begin{equation}
\tilde{H}_\mu (x;q,t)=
\sum_{\lambda\vdash |\mu|}
\left(
\sum_{p\in\mathcal{P}(\lambda )}
q^{\mathrm{inv}_\mu(p)}
t^{\mathrm{maj}_\mu(p)}
\right)
m_{\lambda}(x),
\end{equation}
where $\lambda$ runs over all partitions of size $|\mu|$.
\end{corollary}

To find combinatorial interpretation of the
Kostka--Macdonald polynomials
$\tilde{K}_{\lambda,\mu}(q,t)$
remains significant open problem.
Among many important partial results
about this problem, we would like
to mention the following theorem
also due to Haglund--Haiman--Loehr:
\begin{theorem}[\cite{HHL} Proposition 9.2]
If $\mu_1\leq 2$, we have
\begin{equation}\label{eq:HHL_kosmac}
\tilde{K}_{\lambda,\mu}(q,t)=
\sum_{p\in\mathcal{P}_+(\lambda )}
q^{\mathrm{inv}_\mu (p)}
t^{\mathrm{maj}_\mu (p)},
\end{equation}
where $\mathcal{P}_+(\lambda )$ is the set of all highest weight
elements of $\mathcal{P}(\lambda )$ according to the reading order
explained in Eq.(\ref{eq:Haglundreading}).
\end{theorem}
It is interesting to compare this formula with the
formula obtained by S. Fishel \cite{Fi},
see also \cite{Kir1},~\cite{KiSh}.

Concerning validity of the formula
Eq.(\ref{eq:HHL_kosmac}),
we state the following conjecture.
\begin{conjecture}\label{conj:HHL_kostka}
Explicit formula for
the Kostka--Macdonald polynomials
\begin{equation}
\tilde{K}_{\lambda,\mu}(q,t)=
\sum_{p\in\mathcal{P}_+(\lambda )}
q^{\mathrm{inv}_\mu (p)}
t^{\mathrm{maj}_\mu (p)}.
\end{equation}
is valid if and only if at least
one of the following two conditions
is satisfied.
\begin{enumerate}
\item[(i)]
$\mu_1\leq 3$ and $\mu_2\leq 2$.
\item[(ii)]
$\lambda$ is a hook shape.
\end{enumerate}
\end{conjecture}

\section{Generating function of tau functions}\label{sec:gentau}
\label{sec:sumoftau}
In \cite{KiSa}, we give an elementary proof for special case
$t=1$ of the formula Eq.(\ref{eq:HHL})
in the following form.

\begin{theorem}\label{th:main}
Let $\alpha$ be a composition and $\mu$ be a partition of
the same size. Then,
\begin{equation}\label{eq:main1}
\sum_{p\in\mathcal{P}(\alpha)}
q^{\mathrm{maj}_{\mu'} (p)}=
\sum_{\eta\vdash |\mu|}K_{\eta,\alpha}~K_{\eta,\mu}(q,1).
\end{equation}
\end{theorem}

\begin{conjecture}\label{con:main}
Let $\alpha$ be a composition and $\mu$ be a partition of
the same size. Then,
\begin{equation}\label{eq:main}
q^{-\sum_{i> r}\alpha_i}
\sum_{p\in\mathcal{P}(\alpha)}
q^{\tau_\mu^{r,1} (p)}=
\sum_{\eta\vdash |\mu|}K_{\eta,\alpha}\tilde{K}_{\eta,\mu}(q,1).
\end{equation}
\end{conjecture}
This conjecture contains Conjecture 5.8 of \cite{KiSa} and
Theorem \ref{th:main} above
as special cases $r=1$ and $r=\infty$, respectively.
Also, extensions for paths of more general representations
without partition $\mu$ are discussed in Section 5.3 of \cite{KiSa}.

\begin{example}
Let us consider case $\alpha =(4,1,1)$
and $\mu =(4,2)$.
The following is a list of paths $p$ and the
corresponding value of tau function
$\tau_{(4,2)}^{2,1}(p)$.
For example, the top left corner
\fbox{$111123$\hspace{3mm}$1$} means
$p=\fbox{1}\otimes\fbox{1}\otimes\fbox{1}\otimes
\fbox{1}\otimes\fbox{2}\otimes\fbox{3}$
and $\tau_{(4,2)}^{2,1}(p)=1$.
$$\begin{array}{|cc|cc|cc|cc|cc|cc|}
\hline
111123& 1&
111132& 2&
111213& 2&
111231& 3&
111312& 2&
111321& 1\\\hline
112113& 3&
112131& 4&
112311& 3&
113112& 3&
113121& 2&
113211& 2\\\hline
121113& 4&
121131& 5&
121311& 4&
123111& 5&
131112& 4&
131121& 3\\\hline
131211& 4&
132111& 3&
211113& 1&
211131& 2&
211311& 1&
213111& 2\\\hline
231111& 3&
311112& 5&
311121& 4&
311211& 5&
312111& 6&
321111& 4\\\hline
\end{array}$$
Summing up, LHS of Eq.(\ref{eq:main}) is
$$q^{-1}\sum_{p\in\mathcal{P}((4,1,1))}
q^{\tau_{(4,2)}^{2,1} (p)}=
q^5+4 q^4+7 q^3+7 q^2+7 q+4$$
which coincides with the RHS of Eq.(\ref{eq:main}).
Compare this with $\tau^{1,1}_{(4,2)}$ data for the same set of paths
at Example 5.9 of \cite{KiSa}.
\end{example}

\vspace{5mm}
\noindent
{\bf Acknowledgements:}
The work of RS is supported by Grant-in-Aid
for Scientific Research (No.21740114), JSPS.

\end{document}